\input amssym
\def\bf{\fam\bffam\tenbf}
 
\hsize=15,4truecm
\vsize=22.5truecm
\advance\voffset by 1truecm
\mathsurround=1pt

\def\chapter#1{\par\bigbreak \centerline{\bf #1}\medskip}

\def\section#1{\par\bigbreak {\bf #1}\nobreak\enspace}

\def\sqr#1#2{{\vcenter{\hrule height.#2pt
      \hbox{\vrule width.#2pt height#1pt \kern#1pt
         \vrule width.#2pt}
       \hrule height.#2pt}}}

\def\k{\kappa}
\def\o{\omega}

\def\d{\delta}

\def\l{\lambda}

\def\a{\alpha}

\def\g{\gamma}

\def\n{\eta}


\def\A{{\cal A}}

\def\M{\hbox{\bf M}}


\def\th #1 #2. #3\par\par{\medbreak{\bf#1 #2.
\enspace}{\sl#3\par}\par\medbreak}
\def\rem #1 #2. #3\par{\medbreak{\bf #1 #2.
\enspace}{#3}\par\medbreak}
\def\proof{{\bf Proof}.\enspace}
\def\sqr#1#2{{\vcenter{\hrule height.#2pt
      \hbox{\vrule width.#2pt height#1pt \kern#1pt
         \vrule width.#2pt}
       \hrule height.#2pt}}}
\def\eop{\mathchoice\sqr34\sqr34\sqr{2.1}3\sqr{1.5}3}

                                                                     %
                                                                     %
\newdimen\refindent\newdimen\plusindent                              %
\newdimen\refskip\newdimen\tempindent                                %
\newdimen\extraindent                                                %
                                                                     %
                                                                     %
\def\ref#1 #2\par{\setbox0=\hbox{#1}\refindent=\wd0                  %
\plusindent=\refskip                                                 %
\extraindent=\refskip                                                %
\advance\extraindent by 30pt                                         %
\advance\plusindent by -\refindent\tempindent=\parindent 
\parindent=0pt\par\hangindent\extraindent 
{#1\hskip\plusindent #2}\parindent=\tempindent}                      %
\refskip=\parindent                                                  %
                                                                     %

\def\empty{\emptyset}

\def\raj{\restriction}

\def\da{\downarrow}

\def\nda{\mathrel{\lower0pt\hbox to 3pt{\kern3pt$\not$\hss}\downarrow}}
\def\nDa{\mathrel{\lower0pt\hbox to 3pt{\kern3pt$\not$\hss}\Downarrow}}
\def\nbot{\mathrel{\lower0pt\hbox to 4pt{\kern3pt$\not$\hss}\bot}}
\def\ekom{\mathrel{\lower3pt\hbox to 0pt{\kern3pt$\sim$\hss}\mapsto}}

\def\anR{\mathrel{\lower1pt\hbox to 2pt{\kern3pt$R$\hss}\not}}
\def\anoR{\mathrel{\lower1pt\hbox to 2pt{\kern3pt$\overline{R}$\hss}\not}}

\def\anRm{\mathrel{\lower1pt\hbox to 2pt{\kern3pt$R^{-}$\hss}\not}}

\def\ndda{\mathrel{\lower0pt\hbox to 1pt{\kern3pt$\not$\hss}\downdownarrows}}

\null
\vskip 2truecm
\centerline{\bf ON THE NUMBER OF ELEMENTARY SUBMODELS OF}
\centerline{\bf AN UNSUPERSTABLE HOMOGENEOUS STRUCTURE}
\vskip 1truecm
\centerline{Tapani Hyttinen and Saharon Shelah$^{*}$}
\vskip 2.5truecm

\chapter{Abstract}
\bigskip

We show that if $\M$ is a stable unsuperstable homogeneous
structure, then for most $\k <\vert\M\vert$,
the number of elementary submodels of $\M$
of power $\k$ is $2^{\k}$.

\vskip 2.5truecm

Through out this paper we assume that $\M$ is a stable unsuperstable
homogeneous model such that $\vert\M\vert$ is strongly
inaccessible (= regular and strong limit).
We can drop this last assumption if instead of all elementary submodels of
$\M$ we study only suitably small ones. Notice also
that we do not assume that $Th(\M )$ is stable.
We assume that the reader is
familiar with [HS] and use all the notions and results
of it freely. In [Hy1] a strong nonstructure theorem was proved
for the elementary submodels of $\M$ assuming the existence of
Skolem-functions. In this paper we drop the assumption on the
Skolem-functions and prove the following nonstructure theorem.

\th 1 Theorem. Let $\l$ be the least regular cardinal $\ge\l (\M )$.
Assume $\k$ is an uncountable regular cardinal
($<\vert\M\vert$) such that $\k >\l$ and $\k^{\o}=\k$.
Then there are models (=elementary submodels of $\M$) $\A_{i}$, $i<2^{\k}$,
such that for all $i<2^{\k}$, $\vert\A_{i}\vert =\k$ and
for all $i<j<2^{\k}$, $\A_{i}\not\cong\A_{j}$.

See [Hy1] for nonstructure results in the case $\M$ is unstable.

We prove Theorem 1 in a serie of lemmas.
Let $\l$ and $\k$ be as in Theorem 1.
By
$\l$-saturated, $\l$-primary etc., we mean $F^{\M}_{\l}$-saturated,
$F^{\M}_{\l}$-primary etc. Notice that $\M$ is $\l$-stable.

\vskip 1truecm

\noindent
$*$ Research supported by the United States-Israel Binational
Science Foundation. Publ. 632.

\vfill
\eject

The notion $\l$-construction
(=$F^{\M}_{\l}$-construction) is defined as general $F$-con\-struc\-tion
is defined in [Sh].

\th 2 Lemma. Assume $(C,\{ a_{i}\vert\ i<\a\} ,\{ A_{i}\vert\ i<\a\} )$
is a $\l$-construction and $\sigma$ is a
permutation of $\a$. Let $b_{i}=a_{\sigma (i)}$ and
$B_{i}=B_{\sigma (i)}$. If for all $i<\a$,
$B_{i}\subseteq C\cup\{ b_{j}\vert\ j<i\}$, then
$(C,\{ b_{i}\vert\ i<\a\} ,\{ B_{i}\vert\ i<\a\} )$
is a $\l$-construction.

\proof Exactly as [Sh] IV Theorem 3.3. $\eop$

We write $\k^{\le\o}$ for
$\{ \n :\a\rightarrow\k\vert\ \a\le\o\}$,
$\k^{<\o}$ and $\k^{\o}=\k^{=\o}$ are defined similarly
(of course these have also the other meaning, but it will
be clear from the context, which one we mean).
Let $J\subseteq\k^{\le\o}$ be such that it is closed under initial
segments.
If $\n ,\xi\in J$ then by $r'(\n ,\xi )$
we mean the longest element of $J$ which is an initial
segment of both $\n$ and $\xi$. If $u,v\in I=P_{\o}(J)$ (=the set
of all finite subsets of $J$)
then by $r(u,v)$ we mean the largest set $R$ which satisfies

(i) $R\subseteq\{ r'(\n ,\xi )\vert\ \n\in u,\ \xi\in v\}$

(ii) if $u,v\in R$ and $u$ is an initial seqment of $v$, then $u=v$.

\noindent
We order $I$
by $u\le v$ if for every $\n\in u$ there is $\xi\in v$ such that
$\n$ is an initial seqment of $\xi$ i.e.
$r(u,v)=r(u,u)$ ($=\{\n\in u\vert\ \neg\exists\xi\in u
(\n\ \hbox{\rm is a proper initial segment of}\ \xi)\}$).

\th 3 Definition. Assume $J\subseteq\k^{\le\o}$
is closed under initial segments
and $I=P_{\o}(J)$. Let $\Sigma=\{ A_{u}\vert\ u\in I\}$
be an indexed family of subsets of $\M$ of power $<\vert\M\vert$.
We say that $\Sigma$ is strongly independent if

(i) for all $u,v\in I$, $u\le v$ implies $A_{u}\subseteq A_{v}$,

(ii) if $u,u_{i}\in I$, $i<n$, and $B\subseteq\cup_{i<n}A_{u_{i}}$
has power $<\l$, then there is an automorphism
$f=f^{\Sigma,B}_{(u,u_{0},...,u_{n-1})}$ of $\M$ such
that $f\raj (B\cap A_{u})=id_{B\cap A_{u}}$
and $f(B\cap u_{i})\subseteq A_{r(u,u_{i})}$.

The model construction in Lemma 4 belowe is a generalized version
of the construction used in [Sh] XII.4.

\th 4 Lemma. Assume that $\Sigma=\{ A_{u}\vert\ u\in I\}$,
$I=P_{\o}(J)$,
is strongly independent. Then there are sets $\A_{u}\subseteq\M$,
$u\in I$, such
that

(i) for all $u,v\in I$, $u\le v$ implies $\A_{u}\subseteq\A_{v}$,

(ii) for all $u\in I$, $\A_{u}$ is $\l$-primary over $A_{u}$,
(and so by (i),
$\cup_{u\in I}\A_{u}$ is a model),

(iii) if $v\le u$, then
$\A_{u}$ is $\l$-atomic (=$F^{\M}_{\l}$-atomic)
over $\cup_{u\in I}A_{u}$ and
$\l$-primary over
$\A_{v}\cup A_{u}$,

(iv) if $J'\subseteq J$ is closed under initial segments
and $u\in P_{\o}(J')$, then
$\cup_{v\in P_{\o}(J')}\A_{v}$ is $\l$-constructible over
$\A_{u}\cup\bigcup_{v\in P_{\o}(J')}A_{v}$.

\proof Let $\{ u_{i}\vert\ i<\a^{*}\}$ be an enumeration of
$I$ such that $u\le v$ and $v\not\le u$ implies
$i<j$. It is easy to see that we can choose $\a$,
$\g_{i}<\a$ for $i<\a^{*}$, $a_{\g}$ and
$B_{\g}$ for $\g<\a$, and $s:\a\rightarrow I$ so that

(a) $\g_{0}=0$ and
$(\g_{i})_{i<\a^{*}}$ is increasing and continuous,

(b) if $\g_{i}\le\g <\g_{i+1}$, then $s(\g )=u_{i}$,

(c) for all $\g <\a$, $\vert B_{\g}\vert <\l$ and
if we write for $\g\le\a$, $A^{\g}_{u}=A_{u}\cup
\{ a_{\d}\vert\ \d <\g ,\ s(\d )\le u\}$,
then $B_{\g}\subseteq A^{\g}_{s(\g )}$,

(d) for all $\g <\a$,
if we write $A^{\g}=\cup_{u\in I}A^{\g}_{u}$, then
$t(a_{\g},B_{\g})$ $\l$-isolates $t(a_{\g},A^{\g})$,

(e) for all $i<\a^{*}$,
there are no $a$ and $B\subseteq A^{\g_{i+1}}_{u_{i}}$
of power $<\l$ such that $t(a,B)$ $\l$-isolates
$t(a,A^{\g_{i+1}})$,

(f) if $a_{\d}\in B_{\g}$, then $B_{\d}\subseteq B_{\g}$.

\noindent
For all $u\in I$, we define $\A_{u}=A^{\a}_{u}$.
We show that these
are as wanted.

(i) follows immediately from the definitions and for (ii)
it is enough to prove the following
claim (Claim (III) implies (ii) easily).

{\bf Claim.} For all $i<\a^{*}$,

(I) $\Sigma_{i}=\{ A^{\g_{i}}_{u}\vert\ u\in I\}$
is strongly independent, we write
$f^{i,B}_{(u,u_{0},...,u_{n-1})}$ instead of
$f^{\Sigma_{i},B}_{(u,u_{0},...,u_{n-1})}$,

(II) the functions $f^{i,B}_{(u,u_{0},...,u_{n-1})}$ can be chosen
so that
if $j<i$, $u,u_{k}\in I$, $k<n$, $B\subseteq\cup_{i<n}A^{\g_{i}}_{u_{k}}$
has power $<\l$ and $a_{\g}\in B$ implies $B_{\g}\subseteq B$
and $B'=B\cap A^{\g_{j}}$, then
$f^{i,B}_{(u,u_{0},...,u_{n-1})}\raj B'=
f^{j,B'}_{(u,u_{0},...,u_{n-1})}\raj B'$,

(III) if $j<i$, then $A^{\g_{j+1}}_{u_{j}}$ is $\l$-saturated,

\proof Notice that if $a_{\g}\in A^{\d}_{u}\cap A^{\d}_{v}$,
then $a_{\g}\in A^{\d}_{r(u,v)}$. Similarly we see that
the first half of (I) in the claim is always true
(i.e. if $u\le v$ then for all $\d <\a$, $A^{\d}_{u}\subseteq A^{\d}_{v}$.)
We prove the rest by induction on $i<\a^{*}$.
We notice first that it is enough to prove the existence
of $f^{i,B}_{(u,u_{0},...,u_{n-1})}$ only in the case
when $B$ satisfies

(*) if $a_{\g}\in B$, then $B_{\g}\subseteq B$.

For $i=0$, there is nothing to prove. If $i$ is limit,
then the claim follows easily from the induction assumption
(use (II) in the claim).
So we assume that the claim holds for $i$ and prove it
for $i+1$. We prove first (I) and (II). For this let
$u,u_{k}\in I$, $k<n$, and $B\subseteq\cup_{k<n}A^{\g_{i+1}}_{u_{k}}$
be of power $<\l$ such that (*) above is satisfied.
If for all $k<n$, $s(\g_{i})\not\le u_{k}$, then
(I) and (II) in the claim follow
immediately from the induction assumption. So we may assume that
$s(\g_{i})\le u_{0}$.
Let $B'=B\cap(\cup_{k<n}A^{\g_{i}}_{u_{k}})$. By the induction
assumption there is an automorphism
$f=f^{i,B'}_{(u,u_{0},...,u_{n-1})}$ of $\M$ such that
$f\raj (B'\cap A^{\g_{i}}_{u})=id_{B'\cap A^{\g_{i}}_{u}}$
and $f(B'\cap A^{\g_{i}}_{u_{k}})\subseteq A^{\g_{i}}_{r(u,u_{k})}$.
If $s(\g_{i})\le u$, then, by (*) and (d) in the construction,
we can find an automorphism
$g=f^{i+1,B}_{(u,u_{0},...,u_{n-1})}$ of $\M$ such that
$g\raj B'=f\raj B'$ and $g\raj (B-B')=id_{B-B'}$. Clearly this is as wanted.

So we may assume that $s(\g_{i})\not\le u$. Since $s(\g_{i})\le u_{0}$,
$u_{0}\not\le r(u,u_{0})$.
By the choise of the enumeration of $I$ there is
$j<i$ such that $u_{j}=r(u,u_{0})$. Then by the induction assumption
(part (III)),
$A^{\g_{i+1}}_{u_{j}}=A^{\g_{i}}_{u_{j}}=A^{\g_{j+1}}_{u_{j}}$
is $\l$-saturated and by the choise of $f$,
$f(B'\cap A^{\g_{i}}_{u_{0}})\subseteq A^{\g_{i}}_{u_{j}}$.
So by (d) in the construction and (*) above,
there are no difficulties
in finding the required automorphism $f^{i+1,B}_{(u,u_{0},...,u_{n-1})}$.

So we need to prove (III): For this it is enough to show that
$A^{\g_{i+1}}_{u_{i}}$ is $\l$-saturated. Assume not.
Then there are $a$ and $B$ such that
$B\subseteq A^{\g_{i+1}}_{u_{i}}$,
$\vert B\vert <\l$
and $t(a,B)$ is not realized in $A^{\g_{i+1}}_{u_{i}}$.
Since $\l\ge\l (\M )$, there are $b$ and $C$ such that
$B\subseteq C\subseteq A^{\g_{i+1}}_{u_{i}}$,
$\vert C\vert <\l$, $t(b,B)=t(a,B)$ and
$t(b,C)$ $\l$-isolates $t(b,A^{\g_{i+1}}_{u_{i}})$.
But since (I) in the claim holds for $i+1$,
$t(b,C)$ $\l$-isolates $t(b,A^{\g_{i+1}})$.
This contradicts (e) in the construction. $\eop$ Claim

(iii) and (iv) follow immediately from the construction,
Claim (III) and Lemma 2.
$\eop$

Since $\M$ is unsuperstable, by [HS] Lemma 5.1, there are
$a$ and $\l (\M )$-saturated models $\A_{i}$, $i<\o$,
such that

(i) if $j<i<\o$, then $\A_{j}\subseteq\A_{i}$,

(ii) for all $i<\o$, $a\nda_{\A_{i}}\A_{i+1}$.

\noindent
It is easy to see that we may choose the models $\A_{i}$ so
that they are $\l$-saturated and of power $\l$.
Let $\A_{\o}$ be $\l$-primary over $a\cup\bigcup_{i<\o}\A_{i}$.
As in [Hy1] Chapter 1, for all $\n\in\k^{\le\o}$,
we can find $\A_{\n}$ such that

(a) for all $\n\in\k^{\le\o}$, there is an automorphism $f_{\n}$ of $\M$
such that $f_{\n}(\A_{length(\n )})$ $=\A_{\n}$,

(b) if $\n$ is an initial segment of $\xi$, then
$f_{\xi}\raj\A_{length(\n )}=f_{\n}\raj\A_{length(\n )}$,

(c) if $\n\in\k^{<\o}$, $\a\in\k$ and
$X$ is the set of those $\xi\in\k^{\le\o}$ such that
$\n\frown(\a )$ is an initial segment of $\xi$, then
$$\cup_{\xi\in X}\A_{\xi}\da_{\A_{\n}}\cup_{\xi\in(\k^{\le\o}-X)}\A_{\xi}.$$

\noindent
For all $\n\in\k^{\o}$, we let $a_{\n}=f_{\n}(a)$.

\th 5 Lemma. Assume $\n\in\k^{<\o}$, $\a\in\k$ and
$X$ is the set of those $\xi\in\k^{<\o}$ such that
$\n\frown(\a )$ is an initial segment of $\xi$. Let
$B\subseteq\cup_{\xi\in(\k^{\le\o}-X)}\A_{\xi}$ and
$C\subseteq\cup_{\xi\in X}\A_{\xi}$ be of power $<\l$.
Then there is $C'\subseteq\A_{\n}$ such that
$t(C',B)=t(C,B)$.

\proof By [Hy2] Lemma 8 (or [HS] Lemma 3.15 plus little work)
we can find $D\subseteq\A_{\n}$ of power $<\l$ such that
for all $b\in B$, $t(b,\A_{\n}\cup C)$ does not split over $D$.
So if we choose $C'\subseteq\A_{\n}$ so that
$t(C',D)=t(C,D)$, then $C'$ is as wanted. $\eop$

\th 6 Lemma. Assume $J\subseteq\k^{\le\o}$ and $I=P_{\o}(J)$. For all
$u\in I$, define $A_{u}=\cup_{\n\in u}\A_{\n}$. Then
$\{ A_{u}\vert\ u\in I\}$ is strongly independent.

\proof Follows immediately from Lemma 5. $\eop$

Let $S\subseteq\{\a <\k\vert\ cf(\a )=\o\}$.
By $J_{S}$ we mean the set
$$\k^{<\o}\cup\{ \n\in\k^{\o}
\vert\ \n\ \hbox{\rm is strictly increasing and}
\ \cup_{i<\o}\n (i)\in S\} .$$
Let $I_{S}=P_{\o}(J_{S})$ and $\A_{S}$ be the model given
by Lemmas 4 and 6 for $\{ A_{u}\vert\ u\in I_{S}\}$.

\th 7 Lemma.

(i) Assume $\n\in \k^{<\o}$, $u\in I_{S}$, $\a <\k$,
$\{\n\}\le u$ and
$\{\n\frown (\a )\}\not\le u$. Let $X$ be the set of theose
$\xi\in J_{S}$ such that $\n\frown (\a )$ is
an initial segment of $\xi$. Then
$$\cup_{\xi\in X}\A_{\xi}\da_{\A_{u}}\cup_{\xi\in J_{S}-X}\A_{\xi}.$$

(ii) Assume $\a\in\k$, $u\in I_{S}$
and $v\in P_{\o}(J_{S}\cap\a^{\le\o})$ is maximal
such that $v\le u$. Then
$$\A_{u}\da_{\A_{v}}\cup_{w\in P_{\o}(J_{S}\cap\a^{\le\o})}\A_{w}.$$

\proof (i): Let $C=\cup_{\xi\in X}\A_{\xi}$. By (c) in the
definition of $\A_{\xi}$, $\xi\in\k^{\le\o}$,
there is $C'$ such that
$t(C',\cup_{\xi\in J_{S}-X}\A_{\xi})=t(C,\cup_{\xi\in J_{S}-X}\A_{\xi})$
and $C'\da_{\A_{\n}}\A_{u}\cup\bigcup_{\xi\in J_{S}-X}\A_{\xi}$.
So the claim follows from the first half of Lemma 4 (iii).

(ii): By (i),
$A_{u}\da_{\A_{v}}\cup_{w\in P_{\o}(J_{S}\cap\a^{\le\o})}A_{w}$
from which the claim follows by Lemma 4 (iii) and (iv). $\eop$

\th 8 Lemma. Assume $S,R\subseteq\{ \a <\k\vert\ cf(\a )=\o\}$
are such that $(S-R)\cup (R-S)$ is stationary.
Then $\A_{S}$ is not isomorphic to $\A_{R}$.

\proof Assume not. Let $f:\A_{S}\rightarrow\A_{R}$
be an isomorphism. We write $I^{\a}_{S}$ for the set of
those $u\in I_{S}$, which satisfy that for all
$\xi\in u$, $\cup_{i<length(\xi )}\xi (i) <\a$.
$I^{\a}_{R}$ is defined similarly.
Then we can find $\a$ and $a_{i}$, $i<\o$, such that
$\n =(a_{i})_{i<\o}$ is strictly increasing, for all $i<\o$,
$f(\cup_{u\in I^{\a_{i}}_{S}}\A_{u})=\cup_{u\in I^{\a_{i}}_{R}}\A_{u}$
and $\a=\cup_{i<\o}\a_{i}\in (S-R)\cup (R-S)$.
Without loss of generality we may assume that $\a\in S-R$,
and so $\n\in J_{S}-J_{R}$.
Let $\A^{\a_{i}}_{S}=\cup_{u\in I^{\a_{i}}_{S}}\A_{u}$ and
$\A^{\a_{i}}_{R}=\cup_{u\in I^{\a_{i}}_{R}}\A_{u}$.
Then it easy to see that for all $i<\o$,
$a_{\n}\nda_{\A^{\a_{i}}_{S}}\A^{\a_{i+1}}_{S}$
(use [HS] Lemma 3.8 (iii)).
So there is $u\in I_{R}$ such that for all $i<\o$,
$\A_{u}\nda_{\A^{\a_{i}}_{R}}\A^{\a_{i+1}}_{R}$.
Since $\a\not\in R$, this contradicts Lemma 7 (ii). $\eop$

We can now prove Theorem 1: By [Sh] Appendix 1 Theorem 1.3 (2)
and (3), there are stationary
$S_{i}\subseteq\{\a <\k\vert\ cf(\a )=\o\}$,
$i<\k$, such that for all $i<j<\k$, $S_{i}\cap S_{j}=\empty$.
For all $X\subseteq\k$, let $\A_{X}=\A_{\cup_{i\in X}S_{i}}$.
Then by Lemma 8, if $X\ne X'$, then $\A_{X}$ is not
isomorphic to $\A_{X'}$. Since $\k^{\o}=\k$,
$\vert\A_{X}\vert =\k$. $\eop$ Theorem 1.

\chapter{References.}

\item{[Hy1]} T. Hyttinen, On nonstructure of elementary submodels
of an unsuperstable homogeneous structure,
Mathematical Logic Quarterly, to appear.

\item{[Hy2]} T. Hyttinen, Generalizing Morley's theorem, to appear.

\item{[HS]} T. Hyttinen and S. Shelah, Strong splitting in stable
homogeneous models, to appear.

\item{[Sh]} S. Shelah, Classification Theory, Stud. Logic Found. Math.
92, North-Holland, Amsterdam, 2nd rev. ed., 1990.

\bigskip

Tapani Hyttinen

Department of Mathematics

P.O. Box 4

00014 University of Helsinki

Finland

\medskip

Saharon Shelah

Institute of Mathematics

The Hebrew University

Jerusalem

Israel

\medskip

Rutgers University

Hill Ctr-Bush

New Brunswick

New Jersey 08903

U.S.A.

\end